\input amstex
\documentstyle{amsppt}
\magnification=1200
\NoRunningHeads
\NoBlackBoxes
\topmatter 
\title An explicit formula for the determinant of the Abelian integral 
matrix\endtitle
\author A.A.Glutsyuk\endauthor 
\abstract We consider a polynomial $h(x,y)$ in two complex variables 
of degree $n+1\geq2$ with a 
generic higher homogeneous part. The rank of the first homology 
group of its nonsingular level curve $h(x,y)=t$ is $n^2$. To each 1- form 
in the variable plane and a generator of the homology group one associates  
the (Abelian) integral of the form along the generator. The Abelian integral 
is a multivalued function in $t$. For a fixed canonic tuple of 
$n^2$ monomial 1- forms we consider the multivalued $n^2\times n^2$ 
matrix function in $t$ whose elements are the Abelian integrals of the forms 
along the generators. Its determinant does not depend on the choice 
of the generators in the homology group (up to change of sign, which 
corresponds to change of generator system that reverses orientation). 

In 1999 Yu.S.Ilyashenko proved \cite{1} that the determinant of the 
Abelian integral matrix is a polynomial in $t$ of degree $n^2$ 
whose zeroes are the critical values of $h$. 
We give an explicit formula for the determinant. 
\endabstract
\thanks Research supported by CRDF grant RM1-229, by INTAS grant 93-0570-ext, 
by Russian Foundation for Basic Research (RFBR) grant 98-01-00455, by State 
Scientific Fellowship of Russian Academy of Sciences for young scientists, 
by the European Post-Doctoral Institute joint fellowship of Max-Planck Institut f\"ur 
Mathematik (Bonn) and IH\'ES (Bures-sur-Yvette, France)\endthanks
\address Independent University of Moscow, Bolshoi Vlasievskii Pereulok, 11, 
121002 Moscow, Russia\endaddress
\address Steklov Mathematical Institute, Moscow\endaddress
\address Moscow State University, Department of Mathematics\endaddress
\address Present address: Institut des Hautes \'Etudes Scientifiques, 
 35 Route de Chartres, 91440 Bures-sur-Yvette, France\endaddress
\endtopmatter 

\document
\define\om{\omega}
\define\oj{\om_j}
\redefine\c{\Bbb C}
\define\cc{\c^2}
\redefine\a{\alpha}
\define\djk{\det(I_{j,r})}
\define\dnk{\det E_{n,k}}

\head 1. Introduction. Statement of result\endhead

The infinitesimal Hilbert 16-th problem is the following: {\it what may be 
said about the number and the location of limit cycles of a 
two-dimensional polynomial vector field close to a Hamiltonian one?} 

Yu.S.Ilyashenko \cite{1} have partially investigated the following question 
closely related to the infinitesimal Hilbert problem: to find an explicit 
form of the Picard-Fuchs equation for Abelian integrals.

In the present paper for a hamiltonian of degree $n+1\geq2$ and 
$n^2$ canonic monomial 1- forms (2) we calculate explicitly the determinant of the 
corresponding Abelian integral matrix (3) (Theorem 1).   

To state the main result, let us introduce some notations. 

Let $(x,y)$ be coordinates in $\c^2$. 
Let $n\in\Bbb N$, $n\geq1$, $h(x,y)$   be a polynomial of the degree $n+1$, 
$$H(x,y)=\sum_{s=0}^{n+1}h_sx^{n+1-s}y^s\tag1$$ 
be its homogeneous part of the degree $n+1$. 

Everywhere below we suppose that the homogeneous polynomial $H$ is   
{\it generic}: this means that it is not divisible by square of linear function, 
i.e., has $n+1$ distinct zero lines. Let $a_1$,\dots,$a_{n^2}$ be the critical 
values of the polynomial $h$, $t\in\c\setminus\{ a_1,\dots,a_{n^2}\}$. 
Then the first 
integer homology group of the nonsingular level curve $h(x,y)=t$ has 
dimension $n^2$. 

Let $\{1,\dots,n^2\}\to\{(l,m)|\ 0\leq m,l\leq n-1\}$ 
be the lexicographic numeration of the integer pairs 
$(l,m)$: $(0,0)$, $(0,1)$, $(0,2$... Define the following polynomials and 
differential 1- forms in $\cc$: 
$$e_j(x,y)=x^{l(j)}y^{m(j)}; \ \ d(j)=\deg e_j=l(j)+m(j); \ \ \oj=ye_j(x,y)dx.\tag2$$
Let $t\neq a_i$, $\a=\{\a_r\}$, $r=1,\dots,n^2$, be a system of generators of the group 
$H_1(\{ h(x,y)=t\},\Bbb Z)$. Define the Abelian integrals 
$$I_{j,r}(t)=\int_{\a_r}\oj, \ \ j,r=1,\dots,n^2.\tag3$$
The system $\a$ is chosen to depend continuously on $t\neq a_i$; it is 
"multivalued": $\a$ transforms to another system of generators after 
going around a critical value. The new system is obtained from the initial one 
by unimodular linear transformation.    
The Abelian integrals (3) are multivalued holomorphic functions in $t$ with 
branching points $a_i$. They 
form an $n^2\times n^2$  matrix  with the indices $j,r$.

We calculate the determinant $\det(I_{j,r})(t)$. 

\remark{Remark 1} The function $\det(I_{j,r})(t)$ is holomorphic in $t$ and 
{\it single-valued}. It does not depend (up to 
multiplication by -1) on the choice of the system $\a$ of generators: if 
two generator systems $\a$ and $\a'$ define the same orientation 
of the homology space, then the corresponding determinants are equal; otherwise 
they are opposite. 
\endremark

Earlier Yu.S.Ilyashenko proved \cite{1} that the function $\djk$ 
depends only on $H$: 
$$\djk(t)=C(H)\prod_{i=1}^{n^2}(t-a_i); \ \ a_i \ \text{are the critical values 
of} \ h,\tag4$$
where $C(H)\neq0$ (provided that $H$ is generic). Below we calculate the 
constant  $C(H)$ from (4) as a function in $H$. 

\remark{Remark 2}  The nongeneric polynomials form an irreducible 
hypersurface in the 
space of all the homogeneous polynomials (1) (denote this hypersurface by $S$).  For a generic 
$H$ the corresponding value $C(H)$ is well-defined up to sign. As it is 
shown below (Theorem 1), $C(H)$ is a double-valued function in $H$ with 
branching at $S$. 
\endremark   

By $\Sigma(H)$ denote the discriminant of $H$, i.e., the polynomial in 
the coefficients of $H$ whose zero set is  
 the hypersurface $S$ with the multiplicity 1. 
 
 \remark{Remark 3} Consider the decomposition 
$$H(x,y)=h_0\prod_{i=0}^n(x-b_iy)$$ 
of $H$ into product of linear factors. The discriminant $\Sigma(H)$ is a 
degree $2n$ homogeneous polynomial in the coefficients (1) of $H$. It is equal to  
$$\Sigma(H)=h_0^{2n}\prod_{0\leq j<i\leq n}(b_i-b_j)^2.\tag5$$
\endremark

Let $k\in N$, $k\leq n-1$. Let us introduce the following four upper (lower) 
triangular $k\times k$- matrix linear functions in $H$: 
$$A_{n,k}=\left(\matrix & (n+1)h_0 & nh_1 & \dots & (n-k+2)h_{k-1} \\
& 0 & (n+1)h_0 & \dots & (n-k+3)h_{k-2} \\
& \dots  & \dots & \dots & \dots \\
& 0 & \dots & \dots & (n+1)h_0\endmatrix\right),$$
$$B_{n,k}=\left(\matrix & h_n & 0 &\dots & 0 \\
& 2h_{n-1} & h_n & \dots & 0 \\
& \dots & \dots &\dots &\dots\\
& kh_{n-k+1} &(k-1)h_{n-k+2} &\dots & h_n\endmatrix\right),$$
$$C_{n,k}=\left(\matrix & h_1 & 2h_2 &\dots & kh_k\\
& 0 & h_1 &\dots & (k-1)h_{k-1}\\
& \dots &\dots &\dots &\dots\\
& 0 & 0 &\dots & h_1\endmatrix\right),$$
$$D_{n,k}=\left(\matrix & (n+1)h_{n+1} & 0 &\dots & 0\\
& nh_n & (n+1)h_{n+1} & \dots &0\\
& \dots &\dots &\dots &\dots\\
&(n-k+2)h_{n-k+2} & (n-k+3)h_{n-k+3} &\dots & (n+1)h_{n+1}\endmatrix\right).$$

Define the following $2k\times 2k$ matrix function: 
$$E_{n,k}=\left(\matrix & A_{n,k} & B_{n,k}\\
& C_{n,k} & D_{n,k}\endmatrix\right).\tag6$$

\remark{Remark 4} Let $n\in\Bbb N$, $n\geq2$, $H$, $e_j$, $d(j)$, 
$\oj$ be as at the beginning of the 
paper. Let $k\in\Bbb N$, $1\leq k\leq n-1$,  $E_{n,k}(H)$ be the matrix function 
from (6). Then $E_{n,k}(H)$ is degenerate, iff there exists a nonzero linear 
combination $\sum_{d(j)=n+k-1}c_je_j$ that belongs to the gradient 
ideal of $H$ in the polynomial algebra (i.e., is equal to a linear combination 
of the partial derivatives of $H$ with the coefficients in homogeneous 
degree $k-1$ polynomials). Indeed, for any $l\leq k$ let us assign the monomial $x^{n+k-l}y^{l-1}$ 
($x^{k-l}y^{n+l-1}$) to the $l$- th 
(respectively, $k+l$- th) column of $E_{n,k}$. Then by definition, for any $j\leq k$ 
the $j$- th ($k+j$- th) line of the matrix $E_{n,k}$ consists of the coefficients 
of the degree $n+k-1$ polynomial $x^{k-j}y^{j-1}\frac{\partial H}{\partial x}$ 
(respectively, $x^{k-j}y^{j-1}\frac{\partial H}{\partial y}$) at the monomials 
different from $e_i$.  The last $2k$ polynomials are linearly independent and 
form a basis in the space of polynomials of the degree $n+k-1$ that belong 
to the gradient ideal of $H$. This 
follows from genericity of $H$. Therefore, $E_{n,k}(H)$ is nondegenerate, iff 
each  
linear combination of these $2k$ polynomials contains a monomial different from $e_i$ 
(or equivalently, no nontrivial linear combination of the monomials $e_i$ 
with $d(i)=n+k-1$ belongs to the gradient ideal of $H$).   
\endremark 

\proclaim{Theorem 1} Let $h$, $H$ be as at the beginning of the paper, 
$I_{j,r}$ be the corresponding Abelian inegrals (3),  
$C(H)$ be the corresponding function from (4), $E_{n,k}(H)$ be the matrix functions 
from (6). Then 
$$C(H)= c_n(\Sigma(H))^{\frac12-n}\prod_{k=1}^{n-1}\dnk(H),\tag7$$
$$c_n=(-1)^{\frac{n(3n-1)}4}\frac{(2\pi)^{\frac{n(n+1)}2}(n+1)^{\frac{n^2+n-4}2}((n+1)!)^n}
{\prod_{m=1}^{n-1}(m+n+1)!}.\tag8$$
\endproclaim

\define\prdk{\prod_{k=1}^{n-1}\dnk(H)}
Theorem 1 is proved in Section 2.

\example{Example 1} Let us check the statement of Theorem 1 in the simplest 
case, when $n=1$, $h(x,y)=H(x,y)=x^2+y^2$. Then $C(H)=\pi$: 
by definition, the value 
$C(H)$ is equal to the integral of the form $ydx$ along the unit circle 
in the real $(x,y)$ plane; this integral is equal to the area of the unit disc. Now 
let us calculate $C(H)$ by using (7) and (8). One has 
$$C(H)=c_1(\Sigma(H))^{-\frac12} \ \text{by (7)},\tag9$$
$$\Sigma(H)=-4,\ \ c_1=2\pi i \ \text{by (5) and (8) respectively}.$$
Substituting $(\Sigma(H))^{\frac12}=2i$ and the value of $c_1$ to (9) 
yields $C(H)=\pi$, as above.
\endexample

\head 2. Calculation of $C(H)$. Proof of Theorem 1\endhead
\subhead 2.1. Scheme of the proof of Theorem 1\endsubhead 

Without loss of generality, {\it everywhere below we suppose that} $h\equiv H$, 
so $a_j=0$, $C(H)=\djk|_{t=1}$. Let us calculate the function $C(H)$.

Firstly we show that $C(H)$ has the type (7) with $c_n$ independent on $H$:   
 
 \proclaim{Lemma 1} For any $n\in\Bbb N$, $n\geq1$, there exists a 
nonzero constant $c_n$ satisfying (7) for any generic $H$.
\endproclaim

Lemma 1 is proved in Subsections 2.2-2.4.

Then in Subsection 2.5 we calculate the constant $c_n$. 

\subhead 2.2. Homogeneity degrees and orders of zeroes. Proof of Lemma 1
\endsubhead

As it is shown below, Lemma 1 is implied by the two following statements. 

\proclaim{Proposition 1} The function C(H) is homogeneous of degree $-n^2$.
\endproclaim

Proposition 1 is proved in Subsection 2.3. 

\proclaim{Lemma 2} Let $E_{n,k}(H)$ be the matrix functions from (6). 
The function $C(H)$ is divisible by each polynomial $\dnk$. The ratio 
$\frac{C(H)}{\prdk}$ is a multivalued holomorphic function with at most 
double branching at the hypersurface $S=\{\Sigma(H)=0\}$. No its branch 
vanishes outside $S$.
\endproclaim

Lemma 2 is proves in Subsection 2.4.

Let us prove the implication of Lemma 1 from Proposition 1 and Lemma 2.
 The function 
$C(H)$ has at most polynomial growth, as $H$  approaches $S$   
(by definition and a theorem of P.Deligne \cite{2}). Recall that the surface $S$ is 
irreducible.  
Therefore, the product of $C(H)$ and appropriate power (that will be referred to, 
as $-s$) of the polynomial $\Sigma(H)$ is a single-valued 
holomorphic function that does not vanish identically in $S$  
(by the singularity reducibility theorem). 
Hence, $C(H)=c_n\prdk(\Sigma(H))^s$ with $c_n$ independent on $H$. 
The degree of each $\dnk$ is equal to $2k$. This together 
with the previous statements on the homogeneity degrees of $\Sigma(H)$ and $C(H)$ 
(Remark 3 and Proposition 1) implies that $s=\frac12-n$. 
This proves Lemma 1 modulo Proposition 1 and Lemma 2.

\subhead 2.3. Proof of Proposition 1\endsubhead 
Let $b\in\c\setminus0$. Let us compare $C(H)$ and $C(bH)$. 
By definition, for any $t\in\c$ the value at $t$ of the 
function $\djk$ corresponding to a polynomial $h=H$ is equal to the value at 
$bt$ of that corresponding to $bh$, i.e.,  
$\djk(t)=C(H)t^{n^2}=C(bH)(bt)^{n^2}$. Therefore, $C(bH)=b^{-n^2}C(H).$ This 
proves Proposition 1.

\subhead 2.4. Proof of Lemma 2\endsubhead

For the proof of Lemma 2 we consider $\djk=(\djk)|_{t=1}$ as a function 
in variable polynomials $e_j$ of degrees $d(j)\geq n$: for each one of these $j$ let us 
choose arbitrary homogeneous polynomials $p_j(x,y)$ of degree $d(j)$ and put 
$$\oj=p_j(x,y)ydx.\tag10$$ 
For the other $j$ define $\oj$ to be the forms from (2). Let $H$ be a generic homogeneous polynomial of degree 
$n+1$ (in the same sense, as at the beginning of the paper). Put 
$p=(p_j|\ d(j)\geq n)$. Consider 
 the value at $t=1$ of the matrix function (3) corresponding to the new $\oj$. 
 We consider its determinant 
$$\djk|_{t=1}=C(H,p)\tag11$$ 
as a function in   
the coefficients of the polynomials $H$ and $p_j$. By definition, 
$C(H,(e_j|\ d(j)\geq n))=C(H)$. 

For the proof of Lemma 2 we construct an 
extension $Pr(H,p)$ of the product $\prdk$ as a function in the coefficients 
of $H$ and $p_j$ with the following properties: 

1) $C(H,p)$ vanishes iff so does $Pr(H,p)$; 

2) for any fixed generic $H$ the gradient of each one of the functions $C(H,p)$, 
$Pr(H,p)$ in the coefficients of all the polynomials $p_j$ does not vanish in a Zariski open 
dense subset of the zero set $Pr=0$.

This will imply that for any fixed generic $H$ the ratio $\frac{C}{Pr}(H,p)$ is 
holomorphic in $p$ and does 
not vanish (in particular, its value $C(H)$ corresponding to $p_j=e_j$ 
is finite and nonzero). 
This will prove Lemma 2. 

In the proof of Lemma 2 we use the following 

\proclaim{Proposition 2} The value $C(H,p)$ of the function from (11) 
vanishes iff there exist a $k$, $1\leq k\leq n-1$, and a linear combination 
$$\sum_{d(j)=n+k-1}c_j
\frac{\partial(yp_j)}{\partial y} \ \ \text{with at least one}\ c_j\neq0\tag12$$ 
that belongs to the gradient ideal of $H$. 
\endproclaim

Proposition 2 was proved by Yu.S.Ilyashenko 
\cite{1} in a particular case, and the proof remains valid in the general case.

\remark{Remark 5} For any $n\geq2$, $1\leq k\leq n-1$, the number of the 
indices $j$ with $d(j)=n+k-1$ is equal to $n-k$. Indeed, it follows from 
(2) that the number of the monimials $e_j$ of degree $n+k-1$ is equal to $n-k$. 
\endremark

Let us construct the function $Pr$. Firstly we define appropriate extension  
of each factor $\dnk$, $1\leq k\leq n-1$, that vanishes exactly iff a linear 
combination (12) belongs to the gradient ideal of $H$. To do this, let us 
introduce the following 

\definition{Definition 1} Let $n\geq1$, $H$ be a generic polynomial of degree 
$n+1$ (in the same sense, as at the beginning of the paper). Let 
$1\leq k\leq n-1$. 
Let $q_j$, $j=1,\dots,n-k$, be $n-k$ homogeneous polynomials of the degree $n+k-1$. For each 
$j=1,\dots,n+k$ define the polynomial $Y_j$ of the degree $n+k-1$ as follows: 
for $j\leq n-k$ put $Y_j=\partial{(yq_j)}{\partial y}$; for $n-k+1\leq j\leq n$ 
put $Y_j=x^{n-j}y^{j-n+k-1}\frac{\partial H}{\partial x}$; for $n+1\leq j\leq 
n+k$ put $Y_j=x^{k-j+n}y^{j-n-1}\frac{\partial H}{\partial y}$.  
Define $A(k,H,q_1,\dots,q_{n-k})$ to be the $(n+k)\times(n+k)$- matrix 
with columns numerated by the indices $s$, $1\leq s\leq n+k$ (being assigned 
with the monomials $x^{n+k-s}y^{s-1}$) whose 
$j$-th line consists of the coefficients of the polynomial $Y_j$. 
\enddefinition 
\define\grad{\operatorname{grad}}

\remark{Remark 6} The matrix function $A$ from the previous Definition vanishes, 
iff a nonzero linear combination (12) belongs to the gradient ideal of $H$.
\endremark

Let $1\leq k\leq n-1$, $A$ be the matrix function from the previous Definition. 
The extension of the function $\dnk$ we are 
looking for is $\det A(k,H, p_j|\ d(j)=n+k-1)$. Its value at $p_j=e_j$ is 
equal to $\dnk(H)$ up to multiplication by constant independent on $H$ (see 
the next Remark). 

The product  
$$Pr(H,p)=\prod_{k=1}^{n-1} \det A(k,H,(p_j|\ d(j)=n+k-1))$$
is the function $Pr$ we are looking for: as it is shown below, it satisfies 
the statements 1) and 2) from the beginning of the Subsection. 
 
\remark{Remark 7} Let $e_j$, $1\leq j\leq n^2$, be the monomials from (2), 
$1\leq k\leq n-1$. The 
number of the indices $j$ with $d(j)=n+k-1$ is equal to $n-k$: denote them 
 $j_1<\dots<j_{n-k}$. Let $q_s=e_{j_s}$,  
$A=A(k,H,q_1,\dots,q_{n-k})$ be the corresponding matrix function from 
the previous Definition. Let $E_{n,k}$ be the 
matrix function from (6). Then $E_{n,k}$ coincides  
with the unique $2k\times2k$- minor of the 
matrix $A$ such that the complementary minor has nonidentically-vanishing 
determinant, and the latter is equal to a nonzero constant that does not 
depend on $H$. In particular, the value of the function 
$\det A(k,H,q_1,\dots,q_{n-k})$  
at $q_s=e_{j_s}$ is equal to $\dnk$ up to multiplication by nonzero 
constant independent on $H$.
\endremark
 
Let $C(H,p)$ be the function from (11). The value of the function $Pr$ at 
$p_j=e_j$ is equal to $\prdk$ (up to 
multiplication by constant) by definition and the previous Remark. The function 
$Pr(H,p)$ vanishes 
iff so does $C(H,p)$. Indeed, $Pr(H,p)=0$ iff there exists 
a $k$, $1\leq k\leq n-1$, such that a nonzero linear combination (12) belongs 
to the gradient ideal of $H$ (by definition and Remark 6). By Proposition 2, 
this is the case iff $C(H,p)$ vanishes. This proves the statement 1) from the beginning 
of the Subsection. Now for the proof of Lemma 2 it suffices to prove the 
statement 2) from the same place. 

Let us prove that the gradient of the function $C$ does not vanish at 
a generic point of its zero set. A generic point $(H,p=(p_j)_{d(j)\geq n})$ 
of its zero set  possesses the following property: 

\item{(13)} $H$ is generic; there exist a unique $k$, $1\leq k\leq n-1$, and a 
unique nontrivial 
linear combination (12) (up to multiplication by constant) such that the 
last combination belongs to the gradient ideal of $H$. 

The set of points satisfying (13) is Zariski open and dense in the zero set of 
$C$. This follows from Proposition 2 and the statement that for each 
$k=1,\dots,n-1$ the codimension of the gradient ideal of $H$ in the space of 
homogeneous polynomials of the degree $n+k-1$ is equal to $n-k$: then 
for a generic point of the set $C(H,p)=0$ the corresponding system of $n-k$ 
polynomials $\frac{\partial(yp_j)}{\partial y}$ has rank $n-k-1$ modulo 
the gradient ideal, which means exactly that (13) holds. Indeed the 
intersection of the gradient ideal with the  space of degree $n+k-1$ 
homogeneous polynomials is equal to $2k$: by definition, it is generated by 
the $2k$ polynomials 
$x^iy^{k-1-i}\frac{\partial H}{\partial x}$, 
$x^iy^{k-1-i}\frac{\partial H}{\partial y}$, which are linearly independent 
(cf. Remark 4). 
 
Let us fix a generic $H$ and show that $\grad C$ (in the variables $p$) does 
not vanish at points satisfying (13).
 Fix a point satisfying (13). By definition, there is an index $l$ such that the corresponding 
coefficient 
$c_l$ from the  linear combination (12) is nonzero. Let us fix such an $l$ and 
consider that $c_l=1$ without loss of generality. 
Then the gradient of the function $C$ along the space of polynomials $p_l$ 
(with fixed $p_j$ corresponding to $j\neq l$) is nonzero. Indeed, let 
$q_l$ be a homogeneous polynomial of the degree $n+k-1$. The derivative of the 
function $C$ 
in $p_l$ in the direction $q_l$ is equal to its value $C(H,(p_j)_{d(j)\geq n,\ j\neq l}, 
q_l)$ 
at $H$ and the polynomials $p_j$ with $j\neq l$ and $q_l$. This value is nonzero 
for a generic $q_l$. This follows from Proposition 2 and the statement that 
 for a generic $q_l$ there is no linear combination (12) (where $p_l$ is changed 
 to $q_l$) that belongs to the 
 gradient ideal of $H$. Indeed by definition, the codimension of the 
 space of polynomials belonging to the gradient ideal in the space of homogeneous 
 polynomials of the degree $n+k-1$ is equal to the number $n-k$ of the indices $j$ 
 with $d(j)=n+k-1$. By (13) and the assumption that $c_l=1$, no linear 
 combination (12) with $c_l=0$ 
 belongs to the gradient ideal. The linear mapping $q\mapsto\frac{\partial
 (yq)}{\partial y}$ of the space of polynomials of a fixed positive degree is 
 an isomorphism. The two last statements imply that for a generic $q_l$ 
 the polynomials $\frac{\partial(yp_j)}{\partial y}$ with $j\neq l$ and  
 $\frac{\partial (yq_l)}{\partial y}$ are linearly independent modulo the gradient 
 ideal. This proves the statement 2) for the function $C$. The proof of the 
 analogous statement on the function $Pr$ repeats the previous one. 
 Lemma 2 is proved. 
  
\subhead 2.5. Calculation of $c_n$\endsubhead
\define\zz{\Bbb Z_{n+1}\oplus\Bbb Z_{n+1}}
\define\var{\varepsilon}
\define\de{\Delta}

Everywhere below we suppose that $H(x,y)=x^{n+1}+y^{n+1}$. 
To find the constant $c_n$, we calculate the values $C(H)$, $\Sigma=\Sigma(H)$ 
and $\dnk(H)$ at the above polynomial $H$. The constant $c_n$ is expressed 
via them by (7).

To sketch the calculations, let 
us introduce the following notations. Put 
$$\var=e^{\frac{2\pi i}{n+1}},\ \ \ \sigma=\prod_{1\leq l<k\leq n+1}(\var^k-\var^l)^2.\tag14$$ 
For any $j=1,\dots,n^2$ put 
$$I_j=\int_0^1x^{l(j)}(1-x^{n+1})^{\frac{m(j)+1}{n+1}}dx,\tag15$$
$$IP=\prod_{j=1}^{n^2}I_j.\tag16$$

In Subsection 2.5.1 we express $C(H)$ via $\sigma$ and $IP$: we show that 
$$C(H)=\sigma^nIP.\tag17$$

In Subsection 2.5.2 we calculate $\sigma$ and $\Sigma$: we show that 
$$\sigma=(-1)^{\frac{n(n-1)}2}(n+1)^{(n+1)},\ \ \Sigma=(-1)^n\sigma=
(-1)^{\frac{n(n+1)}2}(n+1)^{n+1}.\tag18$$

In Subsection 2.5.3 we calculate $IP$: we show that 
$$IP=\frac{(2\pi)^{\frac{n(n+1)}2}(n+1)^{-\frac{n^2+4n+3}2}((n+1)!)^n}
{\prod_{m=1}^{n-1}(m+n+1)!}.\tag19$$

It follows from definition that 
$$\prdk=(n+1)^{n(n-1)}:\tag20$$
the matrices $E_{n,k}$ corresponding to the polynomial $H$ under consideration 
are diagonal with the diagonal elements equal to $n+1$; so, $\dnk=(n+1)^{2k}$, 
which implies (20). 

Now statement (8) of Theorem 1 follows from (17)-(20) and formula (7) 
proved before. This proves Theorem 1.

\subhead 2.5.1. Calculation of $C(H)$. Proof of (17)\endsubhead 

By definition, $C(H)$ is equal to the value 
of the determinant $\djk$ at $t=1$. Let us calculate the latter. 

Let $F=\{ H(x,y)=1\}$. The fiber $F$ admits 
the action of the group $\zz=\{(l,m)|\ \ l,m=0,\dots,n\}$ by 
multiplication by $\var^l$ and $\var^m$ of the coordinates $x$ and $y$ 
respectively. 

We calculate the value $\djk|_{t=1}$ for appropriate base $\a_1,\dots,\a_{n^2}$ in 
$H_1(F,\Bbb Z)$ (defined below) such that each $\a_j$ with $j>1$ is obtained 
from $\a_1$ by the action of the element $(l(j),m(j))\in\zz$. (This basis is 
completely defined by choice of $\a_1$.)   

To define $\a_1$, we us consider the fiber $F$ as a covering over the $x$- axis 
having the branching points with the $x$- coordinates $\var^j$, $j=0,\dots,n$. 
It is the Riemann surface of the multivalued function 
$(1-x^{n+1})^{\frac1{n+1}}$. 

\definition{Definition 2} Let $(x,y)$ be coordinates in complex plane $\Bbb C^2$, 
$F=\{ x^{n+1}+y^{n+1}=1\}$. Consider the radial segments $[0,1]$ and $[0,\var]$ 
of the branching points $1$ and $\var$ respectively in the $x$- axis; the 
former being oriented from 0 to 1, and the latter being oriented from $\var$ to 
0. Their union is an oriented piecewise-linear curve (denote it by $\gamma$). 
Let $\gamma_0$ and $\gamma_1$ be its liftings to the covering $F$ such that 
$\gamma_0$ contains the point $(0,1)$ and $\gamma_1$ is obtained from $\gamma_0$ 
by multiplication of the coordinate $y$ by $\var$. The curves $\gamma_i$, 
$i=0,1$, are oriented from their common origin $(\var,0)$ to their common end 
$(1,0)$. Define $\a_1\in H_1(F,\Bbb Z)$ to be the homology class represented 
by the union of the 
oriented curve $\gamma_0$ and the curve $\gamma_1$ taken with the inverse 
orientation. 
\enddefinition

\proclaim{Proposition 3} Let $F$, $\a_1$ be as in the previous Definition. 
Let $\a_j\in H_1(F,\Bbb Z)$, $j=2,\dots,n^2$, be the homology classes 
obtained from $\a_1$  by the action of the element $(l(j),m(j))\in\zz$. The 
classes $\a_j$, $j=1,\dots,n^2$, generate the homology group. 
\endproclaim

\demo{Proof} Let $\Gamma=\cup_{j=0}^n[0,\var^j]$ be the union of the radial 
segments of the branching points of the fiber $F$  in the $x$- axis. Let 
$\widetilde{\Gamma}\subset F$ be the  preimage of the set $\Gamma$ under the 
projection of $F$ to the $x$- axis. The set 
$\widetilde{\Gamma}$ is a deformation retract of the fiber $F$ by covering 
homotopy theorem (hence, the inclusion $\widetilde{\Gamma}\to F$ is 
a homotopy equivalence). The group $H_1(\widetilde{\Gamma},\Bbb Z)$ is 
generated by $\a_j$ by construction. Hence, this remains valid for the whole 
fiber $F$. This proves Proposition 3. 
\enddemo

Let us calculate the value $\djk|_{t=1}$ in the basis $\a_j$ from Proposition 3. 
To do this, we use the following 

\remark{Remark 8} Let $\oj$ be the forms (2), $\a_j$ be as in Proposition 3, 
$I_{j,r}$ be the corresponding integrals from (3), $(l(j),m(j))$ be the 
lexicographic integer pair sequence from the beginning of the paper. 
For any $j,r=1,\dots,n^2$ 
$$I_{j,r}=\var^{l(r)(l(j)+1)+m(r)(m(j)+1)}I_{j,1}.\tag21$$
\endremark
Formula (21) implies the following 
\proclaim{Corollary 1} Let $\oj$ be the forms (2), $\a_j$ be as in Proposition 3, 
$(I_{j,r})$ be the 
corresponding matrix of the integrals from (3), $\djk$ be its 
determinant.  Put 
$$I=\prod_{j=1}^{n^2}I_{j,1}.$$
Let $(l(j),m(j))$ be the 
lexicographic integer pair sequence from the beginning of the paper. 
Let $G=(g_{jr})$ be the $n^2\times n^2$- matrix with the elements 
$$g_{jr}=\var^{l(r)(l(j)+1)+m(r)(m(j)+1)}.\tag22$$ 
Then 
$$\djk=I\det G.\tag23$$
\endproclaim

Thus, to calculate $\djk$, it suffices to calculate the expressions 
$I$ and $\det G$ from Corollary 1 (by (23)). They fill be calculated separately. 

Firstly we calculate $\det G$. Below we show that 
$$\det G=(n+1)^{-2n}\sigma^n.\tag24$$ 
 Then we calculate $I$ (at the end of the 
Subsection). We show that 
$$I=(n+1)^{2n}IP.\tag25$$
This will prove (17).

\demo{Calculation of $\det G$. Proof of (24)} 
Define the following $n\times n$ matrix: 
$$Q=(q_{jr})=\left(\matrix & 1 & \var & \dots & \var^{n-1}\\
& 1 & \var^2 & \dots & \var^{2(n-1)}\\
& \dots & \dots & \dots & \dots\\
& 1 & \var^n & \dots & \var^{n(n-1)}\endmatrix\right).$$

To calculate  $\det G$, we use the following periodicity property of 
the matrix $G$: its first $n$ columns are filled in by $n$ copies of the 
matrix $Q$, more precisely,
$$\text{for any}\ \ j,r=1,\dots,n,\ \ s=1,\dots,n-1,\ \ \ \ 
g_{j+sn\ r}=g_{j\ r} =q_{j\ r}.\tag26$$

For any $s=1,\dots,n$ by  $Q_s$ denote the $ns\times ns$ matrix formed by the 
first $ns$ lines and columns of the matrix $G$ 
(thus, $Q_n=G$, $Q_1=Q$). To find $\det G=\det Q_n$, we calculate the 
determinant $\det Q_s$ by induction in $s$. 

We prove the following recurrent formula for $\det Q_s$:
$$\det Q_s=\left(\prod_{l=1}^{s-1}(\var^s-\var^l)\right)^n\det Q\det Q_{s-1}.\tag27$$
\demo{Proof of (27)} For the proof of (27), let us transform 
$Q_s$ by preserving its determinant as follows: for each $j\geq n+1$ we 
subtract the $(j-n)$-th line from the $j$-th one. This yields 
the new matrix (that will be denoted by $Q_s'$)  whose first $n$ columns 
are formed by a single copy of the matrix $Q$ that takes the first $n$ lines 
and zeroes in the other places (by (26)). By definition, 
$$(Q_s')_{j\ r}=g_{j\ r}-g_{j-n\ r}\ \ \text{(we put} \ g_{j\ r}=0\  
\text{for}\ j\leq0).\tag28$$
Let $Q''_s$ be the $n(s-1)\times n(s-1)$- 
matrix obtined from $Q_s'$ by throwing away its first $n$ columns and first $n$ 
lines. By definition, 
$$\det Q_s=\det Q\det Q''_s.\tag29$$ 
Now let us calculate $\det Q''_s$. By definition and (28), 
$$(Q''_s)_{j\ r}=(Q_s')_{j+n, r+n}=g_{j+n\ r+n}-g_{j\ r+n}.$$
 By (22), for any $j,r=1,\dots,n^2-n$
$$g_{j+n\ r}=\var^{l(r)}g_{j\ r}; \ \ g_{j\ r+n}=\var^{l(j)+1}g_{j\ r}.$$ 
This follows from (18) and the relations $l(j+n)=l(j)+1$, $m(j+n)=m(j)$ valid for all $j$ 
(the definition of the sequences $l(j)$ and $m(j)$). Hence, for any $j,r\leq n(s-1)$ 
$$g_{j+n\ r+n}=\var^{(l(r)+1)+(l(j)+1)}g_{j\ r},$$
$$(Q''_s)_{j\ r}=\var^{l(j)+1}(\var^{l(r)+1}-1)g_{j\ r}.$$
In other terms, the matrix $Q''_s$ is obtained from the matrix $Q_{s-1}$ by 
multiplication of its $j$-th column by $\var^{l(j)+1}$ and $r$-th
line by $\var^{l(r)+1}-1$. Therefore, 
$$\det Q''_s=(\prod_{j=1}^{n(s-1)}\var^{l(j)+1})\times(\prod_{r=1}^{n(s-1)}
(\var^{l(r)+1}-1))\det Q_{s-1}=\left((\prod_{l=1}^{s-1}\var^l)\times
(\prod_{k=1}^{s-1}(\var^k-1))\right)^n\det Q_{s-1}.$$
Putting together the product terms corresponding to $l$ and $k$ with $l+k=s$ 
in the previous formula and substituting the latter to (29) yields (27). 
\enddemo

Let us calculate $\det G$. 
Put $\de=\sqrt{\sigma}=\prod_{1\leq l<k\leq n+1}(\var^k-\var^l)$. 
By (27),  
$$\det G=(\det Q)^n\left(\prod_{1\leq l<k\leq n}(\var^k-\var^l)\right)^n.$$
By van der Mond formula,  
$$\det Q=\prod_{1\leq l<k\leq n}(\var^k-\var^l);\ \ \text{so,} \ 
\det G=\left(\prod_{1\leq l<k\leq n}(\var^k-\var^l)\right)^{2n}.\tag30$$ 
The product in (30) is equal to $(n+1)^{-1}\de$ (this statement implies (24)). 
Indeed, by definition, this product is equal to  
$$\de\left(\prod_{1\leq l\leq n}(1-\var^l)\right)^{-1}.$$
The previous statement on (30) follows from the formula
$$\prod_{1\leq l\leq n}(1-\var^l)=n+1:\tag31$$
by definition, its left-hand side is 
the value of the polynomial $\frac{x^{n+1}-1}{x-1}=\sum_{l=0}^nx^l$ at the point 
$x=1$; hence it is equal to $n+1$. Formula (24) is proved.
\enddemo

\demo{Calculation of $I$. Proof of (25)} Let us express $I_{j,1}$ via the 
integral $I_j$ from (15). We  show that 
$$I_{j,1}=(1-\var^{m(j)+1})(1-\var^{l(j)+1})I_j.\tag32$$ 
This together with (31) will imply (25). 

Let $\gamma_0$, $\gamma_1$ be the oriented curves from Definition 2. Then 
$$I_{j,1}=\int_{\a_1}x^{l(j)}y^{m(j)+1}dx=
\int_{\gamma_0}x^{l(j)}y^{m(j)+1}dx-\int_{\gamma_1}x^{l(j)}y^{m(j)+1}dx.\tag33$$ 
 The second integral in the right-hand side of (33) is 
equal to the first one times $\var^{m(j)+1}$ (by definition). Analogously, the first integral 
in its turn is the integral along the segment $[0,1]$ (which  
is equal to $I_j$) minus the one along the segment $[0,\var]$ oriented from 0 
to $\var$. The integral along the last segment is equal to $\var^{l(j)+1}I_j$. 
This together with the two previous statements implies 
(32). Formula (25) is proved. The proof of formula (17) is completed.
\enddemo

\subhead 2.5.2. Calculation of $\sigma$ and $\Sigma$. Proof of (18)\endsubhead
Let us calculate $\sigma$. By definition,
$$\sigma=\left(\prod_{1\leq l<k\leq n+1}(\var^k-\var^l)\right)^2=
(-1)^{\frac{n(n+1)}2}\prod_{1\leq l<k\leq n+1}
((\var^k-\var^l)(\var^l-\var^k))$$ 
$$=(-1)^{\frac{n(n+1)}2}\prod_{1\leq k\leq n+1}\left(\prod_{1\leq l\leq n+1; l\neq k}
(\var^k-\var^l)\right)=(-1)^{\frac{n(n+1)}2}\prod_{1\leq k\leq n+1}\left(\prod_{l=1}^n
\var^k(1-\var^l)\right).$$
Changing the second (inner) product in the right-hand side of the previous formula to 
$(n+1)(\var^k)^n$ (by (31)) yields
$$\sigma=(-1)^{\frac{n(n+1)}2}\var^{\frac{n(n+1)(n+2)}2}(n+1)^{n+1}.$$ 
Substituting $\var^{\frac{n+1}2}=-1$ to the right-hand side of the last formula, 
calculating the resulting power of -1 and changing it to its appropriate 
representative modulo 2 yields the first formula in (18). 

Let us calculate $\Sigma$. By (6), 
$$\Sigma=\left(\prod_{1\leq l<k\leq n+1}\exp(\frac{\pi i}{n+1})
(\var^k-\var^l)\right)^2=
\exp(n(n+1)\frac{\pi i}{n+1})\sigma=(-1)^n\sigma.$$ 
The second formula in (18) is proved.

\subhead 2.5.3. Calculation of $IP$. Proof of (19)\endsubhead
To calculate $IP=\prod_{j=1}^{n^2}I_j$, we firstly express it via appropriate 
values of $B$- and $\Gamma$- functions. 
Recall their definitions:
$$B(a,b)=\int_0^1x^{a-1}(1-x)^{b-1}dx,$$
$$\Gamma(a)=\int_0^{+\infty}x^{a-1}e^{-x}dx.$$
The variable change $u=x^{n+1}$ transforms integral (15) to 
$$\frac1{n+1}\int_0^1u^{\frac{l(j)+1}{n+1}-1}(1-u)^{\frac{m(j)+1}{n+1}}du=
\frac1{n+1}B(\frac{l(j)+1}{n+1},\frac{m(j)+1}{n+1}+1).$$
Therefore,
$$IP=(n+1)^{-n^2}\prod_{0\leq l,m\leq n-1}B(\frac{l+1}{n+1},\frac{m+1}{l+1}+1).
\tag34$$

To calculate the product in the right-hand side of (34), we use the following 
expression of $B$- function via $\Gamma$- function: 
$$B(a,b)=\frac{\Gamma(a)\Gamma(b)}{\Gamma(a+b)}.$$
Therefore, by (34),  
$$IP=(n+1)^{-n^2}\frac{\left(\prod_{l=0}^{n-1}\Gamma(\frac{l+1}{n+1})\right)^n
\left(\prod_{l=0}^{n-1}\Gamma(\frac{l+1}{n+1}+1)\right)
^n}{\prod_{l,m=0}^{n-1}\Gamma(\frac{l+m+2}{n+1}+1)}.\tag35$$
To calculate the products in (35), we use the following 
identities for $\Gamma$- function \cite{3}:
$$\Gamma(n)=(n-1)!\ \ \text{for any}\ \ n\in\Bbb N,$$
$$\prod_{l=0}^n\Gamma(z+\frac l{n+1})=(2\pi)^{\frac n2}(n+1)^{\frac12-(n+1)z}
\Gamma((n+1)z)\ \ \text{(Gauss-Legendre formula)}.\tag36$$

One gets  
$$\prod_{l=0}^{n-1}\Gamma(\frac{l+1}{n+1})=(2\pi)^{\frac n2}(n+1)^{-\frac12},
\tag37$$
$$\prod_{l=0}^{n-1}\Gamma(\frac{l+1}{n+1}+1)=(2\pi)^{\frac n2}(n+1)^{\frac12-(n+2)}
(n+1)!\tag38$$
by applying (36) to $z=\frac1{n+1}$ and $z=\frac{n+2}{n+1}$ respectively and 
subsequent substitutions $\Gamma(1)=1$, $\Gamma(n+2)=(n+1)!$.
Let us calculate the double product in (35). For any fixed $m=0,\dots,n-1$ 
$$\prod_{l=0}^{n-1}\Gamma(\frac{l+m+2}{n+1}+1)=(\Gamma(\frac{m+1}{n+1}+1))^{-1}
(2\pi)^{\frac n2}(n+1)^{\frac12-(m+n+2)}(m+n+1)!$$
by (36) applied to $z=\frac{m+1}{n+1}+1$. Therefore, 
$$\prod_{l,m=0}^{n-1}\Gamma(\frac{l+m+2}{n+1}+1)=\left(\prod_{m=0}^{n-1}
\Gamma(\frac{m+1}{n+1}+1)\right)^{-1}
(2\pi)^{\frac{n^2}2}(n+1)^{\frac n2-\sum_{m=0}^{n-1}(m+n+2)}\prod_{m=0}^{n-1}
(m+n+1)!$$
Substituting formula (38) for the first product in the right-hand side 
of the last formula and summarizing the power of $n+1$ yields 
$$\prod_{l,m=0}^{n-1}\Gamma(\frac{l+m+2}{n+1}+1)=
(2\pi)^{\frac{n^2-n}2}(n+1)^{-\frac{3(n^2-1)}2}
\prod_{m=1}^{n-1}(m+n+1)!.\tag39$$
Substituting (37)-(39) to (35) yields (19). The proof of Theorem 1 is completed.

\head Acknowledgements\endhead

I am grateful to Yu.S.Ilyashenko for attracting my attention to the problem 
and helpful discussions. The paper was written while I was staying at the 
  Institut des Hautes \'Etudes Scientifiques 
(IH\'ES, Bures-sur-Yvette, France). My stay at IH\'ES was 
supported by the European Post-Doctoral Institute fellowship of 
IH\'ES. I wish to thank the Institute for the hospitality and support.

 \head References\endhead
1. Ilyashenko, Yu.S. Restricted infinitesimal Hilbert 16-th problem and 
Picard-Fuchs equations for Abelian integrals. - To appear.

2. Deligne, P.  Equations diff\'erentielles \`a points singuliers r\'eguliers. - 
Lect. Notes Math., V.163. Springer-Verlag, Berlin-New York, 1970.

3. Erdelyi, A. (Editor). Higher transcendental functions, v.1. Bateman 
manuscript project. - McGraw Hill Book Company INC, 1953.

\enddocument